\documentclass[12pt]{article}
\usepackage{amsbsy}

\begin{document}
\def\mapright#1{\smash{\mathop{\longrightarrow}\limits^{#1}} }
\def\mapleft#1{\smash{\mathop{\longleftarrow}\limits^{#1}}}
\def\mapdown#1{\Big\downarrow\rlap{$\vcenter{\hbox{$\scriptstyle#1$}}$}}
\def\mapup#1{\Big\uparrow\rlap{$\vcenter{\hbox{$\scriptstyle#1$}}$}}
\def\mapne#1{\nearrow\rlap{$\vcenter{\hbox{$\scriptstyle#1$} }$}}
\def\mapse#1{\searrow\rlap{$\vcenter{\hbox{$\scriptstyle#1$} }$}}
\def\mapsw#1{\swarrow\rlap{$\vcenter{\hbox{$\scriptstyle#1$} }$}}
\def\mapnw#1{\nwarrow\rlap{$\vcenter{\hbox{$\scriptstyle#1$} }$}}
\def\morph#1{\smash{\mathop{\cong}\limits^{#1}}}
\def\sumup{\displaystyle\sum}
\font\frk=eufm10 at 12pt
\font\frkt=eufm10 at 20pt \def\goth#1{\hbox{\frk #1}}
\def\gotht#1{\hbox{\frkt #1}}
\newtheorem{Proposition}{Proposition}[section]
\newtheorem{Theorem}{Theorem}[section]
\newtheorem{Lemma}{Lemma}[section]
\newtheorem{cory}{Corollary}[section]
\newtheorem{Remark}{Remark}[section]
\newtheorem{define}{Definition}[section]
\baselineskip7mm
\topmargin0in
\date{}

\title{The Lie algebra cohomology of jets\footnote{MSC subject
classification 17B56, 17B65. Written while holding a research
fellowship at the Abdus Salam International Centre for Theoretical
Physics, Trieste, Italy and the Max-Planck Institut f\"ur mathematik,
Bonn, Germany.}}
\author{Yunhyong Kim\\\small Max-Planck Instut fuer Mathematik\\
\small Vivatsgasse 7, 53111 Bonn, Germany\\\small
(e-mail:kimy@mpim-bonn.mpg.de)}
\maketitle

\begin{abstract}
Let $\goth{g}$ be a finite-dimensional complex reductive semi simple Lie
algebra. We present a new calculation of the continuous cohomology of
the Lie algebra $z \goth{g}[[z]]$. In particular, we shall give an
explicit formula for the Laplacian on the Lie algebra cochains, from
which we can deduce that the cohomology in each dimension is a
finite-dimensional representation of $\goth{g}$ which contains any irreducible
representation of $\goth{g}$ at most once.
\begin{description}
\item[Keywords:] Lie algebra, cohomology, semi-infinite forms, affine roots
\end{description}
\end{abstract}

\section{Introduction}
Let $\goth{g}$ be a complex semi-simple Lie algebra.  In this paper, we
shall calculate the cohomology of the Lie algebra $z \goth{g}[[z]]$ of
formal power series (with vanishing constant term) by an infinite
dimensional analog of the method described in the paper by
B. Kostant~\cite{kost}.

The Lie algebras of interest in Kostant's paper are nilpotent Lie
subalgebras $\goth n$ of a finite dimensional semi-simple Lie algebra
$\goth{g}$. Kostant identifies the cohomology of $\goth{n}$ with the
kernel of the Laplace operator on the cochains of $\goth{n}$. The
cochains of $\goth{n}$ can be identified with a summand in the
cochains of the larger algebra $\goth{g}$, and Kostant defines an
operator $\tilde L$ on the cochains of $\goth{g}$ which restricts to
the Laplace operator on the cochains of $\goth{n}$. Calculating the
kernel of $\tilde L$, which turns out to be easier than a direct
calculation of the Laplacian on the Lie algebra cochains of the
nilpotent Lie algebra, yields the cohomology of the Lie subalgebra.

The Lie algebra $\goth{a} = z \goth{g}[[z]]$ with which we are concerned, is
an infinite dimensional topologically nilpotent subalgebra of the algebra
$\goth{g}[[z]][z^{-1}]$ of formal loops in $\goth{g}$. We would like
to emulate Kostant's method in the following way. First we will define
and describe a graded complex of ``semi-infinite forms''. On this
complex, we will define an operator $\tilde L$. The Lie
algebra cochains of $\goth{a}$
will be shown to be a subcomplex of the semi-infinite forms. It
will be proved that the operator $\tilde L$ restricts to the
Laplacian on the subcomplex. We will then give an explicit formula
for $\tilde L$, which will finally enable us to calculate its
kernel and give a description of the cohomology  of $\goth{a}$.

The final result of this paper follows already from the theorem of
H. Garland and J. Lepowsky~\cite{garland}. However, they make use of
the weak Berstein-Gelfand-Gelfand resolution and do not concern
themselves with an explicit description of the Laplacian on the Lie
algebra cochains.  The calculation in this paper gives us an explicit
formula, and is also useful in connection with the smooth cochain
cohomology of loop groups, which I hope to discuss in a separate
paper. The discussion also relates the cohomology to semi-infinite
cohomology which is of independent interest. Finally, I
believe the translation of Kostant's result into an infinite
dimensional setting is appealing in itself, because it illustrates the
power of his method.

To be precise, we will describe the Lie algebra $\goth{a}$ in the
following way: let $G$ be a compact connected simply connected real
Lie group. Let $\goth{g}$ denote the Lie algebra of $G$.  Consider the
Lie algebra $\goth{A}$ consisting of Laurent polynomials of the form

$$\sumup_{p\in {\bf Z}} A_p z^p$$
where $p$ runs over the integers, $A_p$ is in the complexification
$\goth{g}_{\bf C}$ of $\goth{g}$, and such that $A_p =0$ for all but
finite number of $p$'s. Given $A=\sumup_{p\in {\bf Z}} A_p z^p$
and $B= \sumup_{p\in {\bf Z}} B_p z^p$, the Lie bracket is

$$[A, B] = \sumup_{p,q\in {\bf Z}} [A_p, B_q] z^{p+q}.$$
Note that $\goth{A}$ can be decomposed as

$$\goth{A} = \overline{\goth{a}} \oplus \goth{g}_{\bf C}
\oplus \goth{a}$$ 
where $\overline{\goth{a}}$ is the Lie algebra consisting of elements
of the form $\sumup_{k<0} A_k z^k$, and $\goth{a}$ is the Lie algebra
consisting of elements of the form $\sumup_{k >0} A_k z^k$. We would
like to calculate the cohomology of $\goth{a}$ where we consider the
$p$th degree cochains $A^*(\goth{a})$ of the Lie algebra to be complex
multi-linear continuous maps

$$\goth{a}\times \ldots \times \goth{a} \mapright{} {\bf C}.$$

The Lie algebra $\goth{a}$ can be related to the real Lie algebra  quotient $\goth{J}$ of
the Lie algebra $\Omega\goth{g}$ of based loops in $\goth{g}$ by those
whose derivatives vanish to infinite order. The Lie algebra $\goth{J}$ can
be identified with the Lie algebra of jets of loops at the based
point. The cochains on the complexification of $\goth{J}$ is a
subspace of the cochains on $\goth{a}$ and the inclusion induces an
isomorphism on the level of cohomology.  We will come back to all this
with some more detail at the end of this paper.

In Section~\ref{first} we will describe the semi-infinite forms on
$\goth{A}$ and introduce an operator $\tilde{L}$ on these forms. We
will also prove that $\tilde{L}$ restricts to the Laplacian on the
cochains of $\goth{a}$. Then, in section~\ref{second}, we will write
down an explicit formula for $\tilde{L}$, which will enable us to
calculate the kernel of the Laplacian. Section~\ref{third} will
summarise the results following from the formula for the
Laplacian. Finally, in section~\ref{fourth}, we will discuss the
relationship between $\goth{a}$ and the loop group of $G$.

Before, beginning the next section, I would like to make a comment on
notation. There will be many infinite sums in this paper. To avoid
ambiguity, every effort has been made to keep track of indices over
which the sum is to be taken. However, it is to be understood that, as
a rule, repeated indices will be summed overthe integers unless a
restriction has been indicated.

\section{Semi-infinite forms}
\label{first}

We will first define the cochain complex of ``semi-infinite forms'' on
$\goth{A}$. This definition follows the one found
in~\cite{semi-infinite} (Section 1).

Let $c$ be the coxeter number of $\goth{g}$. Let $\langle , \rangle$
be the $\frac{1}{2c}$ times the Killing form. Choose an orthonormal
basis $\{\alpha_i\}_i$ of $\goth{g}_{\bf C}$ with respect to $\langle
, \rangle$. Let $e_{i,k} = \alpha_i z^k$. Then $\{e_{i,k}\}_{i,k}$
form a basis of $\goth{A}$. Denote by $e^{i,k}$ the basis element of
the dual $\goth{A}^*$.  Let us order the basis elements $e^{i,k}$
lexicographically, i.e.

$$e^{i, k} < e^{j, l} \mbox{ if } k<l,\mbox{ or }k = l\mbox{ and
}i< j.$$

For temporary convenience, let us use this ordering to re-index the
basis elements by a single index $\{e_i\}$ and the dual basis $\{e^i\}$.

Define the space $\wedge_\infty^d(\goth{A})$ of semi-infinite forms of
degree $d$ as the space spanned by formal symbols of the form
$\omega = e^{i_1}\wedge e^{i_2}\wedge \ldots \wedge e^{i_p}\wedge
\ldots$ such that there exist $N(\omega)\in {\bf Z}$ for each
$\omega$, so that, for all $p > N(\omega)$, $i_p = p-d$ and such that

$$ e^{i_1}\wedge e^{i_2}\wedge \ldots \wedge e^{i_{p-1}}\wedge
e^{i_{p+1}}\wedge e^{i_p}\wedge e^{i_{p+2}}\wedge \ldots$$
$$= - e^{i_1}\wedge e^{i_2} \wedge \ldots \wedge e^{i_p}\wedge \ldots.$$
Note that elements in $\wedge_\infty^*= \oplus_d \wedge_{\infty}^d$ are
bi-graded: besides the grading by the degree there is a second grading
by the {\it energy} where the energy of $\omega$ is defined as
$\Sigma_p (i_p -p + d)$.

Given any $x =\sumup_i x_i e_i$ in $\goth{A}$ and $x'=\sumup_i
x_ie^i$ in $\goth{A}^*$, we can define operators $\iota(x)$
and $\epsilon(x')$ by

$$\iota(x)(e^{i_1}\wedge \ldots\wedge e^{i_p}\wedge \ldots)$$
$$= \sumup_{p,i} (-1)^{p-1} x_i e^{i_p}(e_i) e^{i_1}\wedge \ldots \wedge{\hat
e^{i_p}}\wedge \ldots$$
where ${\hat e^{i_p}}$ means that the term will be omitted, and

$$\epsilon(\sumup_i x_ie^i)(e^{i_1}\wedge \ldots\wedge e^{i_p}\wedge
\ldots)$$
$$= \sumup_i x_i e^i\wedge e^{i_1}\wedge \ldots \wedge e^{i_p}\wedge
\ldots.$$
In order to simplify the notation, we will write $\iota(e_{i,k})$ and
$\epsilon(e^{i,k})$ as $\iota_{i,k}$ and $\epsilon^{i,k}$
respectively. These operators serve to define
$\wedge_\infty^*(\goth{A})$ as a module of the Clifford algebra on
$\goth{A}\oplus \goth{A}^*$ associated to the pairing $<x, x'>$ for
$x\in \goth{A}$ and $x'\in \goth{A}^*$. That is, the anti-commutator
$[\iota_{i,k}, \epsilon^{j,m}]_+ = \delta_{ij}\delta_{km}$, where
$\delta_{ij} = 0$ if $i\ne j$ and $\delta_{ii} = 1$.

Let

$$:\iota_{i,k}\epsilon^{j,m}: =\left\{\begin{array}{ll}
\iota_{i,k}\epsilon^{j,m} & \mbox{ if }k\leq 0,\\
-\epsilon^{j,m}\iota_{i,k} & \mbox{ if } k > 0 \end{array}\right.$$
For each $x\in \goth{Lg}_{\bf C}$, there is an operator ${\cal L}(x)$  on
$\wedge_\infty^*(\goth{A})$, defined by

$${\cal L}(e_{i,k}) = \sumup C_{iq}^p:\iota_{p,s}\epsilon^{q,s-k}:$$
where $C_{iq}^p$ are the structure constants with respect to the
basis $\{\alpha_p\}$, i.e.,

$$[\alpha_i, \alpha_q] = \sumup_p C_{iq}^p \alpha_p.$$
Since we had taken an orthonormal basis with respect to the killing
form, the structural constants $C_{iq}^p$ are anti-symetric in the
three indexes $i,p,q$. This, along with the identity $[\iota_{i,k},
\epsilon^{j,m}]_+ = \delta_{ij}\delta_{km}$, implies that we could just as
well have written ${\cal L}(e_{i,k}) =\sumup C_{iq}^p
\iota_{p,s}\epsilon^{q,s-k}$.  Again, we will simplify ${\cal
L}(e_{i,k})$ to ${\cal L}_{i,k}$. Note

\begin{Proposition}
The commutator $[\iota_{j.m}, {\cal L}_{i,k}]$ is given by

$$[\iota_{j,m}, {\cal L}_{i,k}] = - \sumup_p C_{ij}^p \iota_{p,m+k}.$$
\end{Proposition}
{\bf Proof:}
 Note that

$$[\iota_{j,m}, {\cal L}_{i,k}] = [\iota_{j,m}, \sumup_{q,p,s} C_{iq}^p
\iota_{p,s} \epsilon^{q,s-k}]$$
$$= \sumup_{p,q,s} C_{iq}^p [\iota_{j,m}, \iota_{p,s}]_+
\epsilon^{q,s-k} - \sumup_{p,q,s\leq 0} C_{iq}^p \iota_{p,s}
[\iota_{j,m},\epsilon^{q,s-k}]_+$$
$$= - \sumup_{p,q,s\leq 0} C_{iq}^p \iota_{p,s} \delta_{jq}\delta_{m(s-k)}$$
$$- \sumup_{p,q,s >0} C_{iq}^p
\delta_{jq}\delta_{m(s-k)}\iota_{p,s}$$
$$= -\sumup_{p} C_{ij}^p \iota_{p,k+m}$$

Likewise, it is just as easy to see that  $[\epsilon^{j,m} , {\cal
L}_{i,k}]= \sumup_q C_{iq}^j \epsilon^{q, m-k}$. Although we have written simply ${\cal L}(e_{i,k}) =\sumup C_{iq}^p \iota_{p,s}\epsilon^{q,s-k}$ in the proof above, note that one must becareful when calculating commutators of infinite sums: given a infite sum $\Sigma_i^\infty a_i$ which converges, its commutator with an element $b$ may not be $\Sigma_i^\infty [a_i, b]$ because the latter may not converge. This is why in Lemma~\ref{hunting} below, we have re-introduced the normal ordering. 

\begin{Proposition}
\label{fish} The operators ${\cal L}_{i,k}$ define a projective
representation of $\goth{A}$ on $\wedge_\infty^*(\goth{A})$.
\end{Proposition}
{\bf Proof of Proposition~\ref{fish}:}
This follows immediately from the following lemma.

\begin{Lemma}
\label{hunting}
The commutator $[{\cal L}_{i,k}, {\cal L}_{j,m}]$ is given by

$$[{\cal L}_{i,k}, {\cal L}_{j,m}] = {\cal L}([e_{i,k}, e_{j,m}]),$$
if $m\ne -k$, and,

$$[{\cal L}_{i,k}, {\cal L}_{j,-k}] = {\cal L}([e_{i,k}, e_{j,-k}])
+ 2c\cdot \delta_{ij} \rangle k.$$

\end{Lemma}
{\bf Proof:}
Assume that $m\geq 0$(if $m\leq 0$, we merely need to
replace $m$ by $-m$).
Note that

$$[{\cal L}_{i,k}, {\cal L}_{j,m}] = [\sumup_{q,p,s} C_{iq}^p
:\iota_{p,s}\epsilon^{q,s-k}:, {\cal L}_{j,m}]$$
$$= - [\sumup_{s >0} C_{iq}^p \epsilon^{q,s-k}\iota_{p,s} ,
{\cal L}_{j,m}] +
[\sumup_{s\leq 0} C_{iq}^p \iota_{p,s}\epsilon^{q,s-k}, {\cal L}_{j,m}]$$
$$= -\sumup_{s >0} C_{iq}^p [\epsilon^{q,s-k}, {\cal L}_{j,m}]\iota_{p,s}
-\sumup_{s>0}C_{iq}^p \epsilon^{q,s-k}[\iota_{p,s}, {\cal L}_{j,m}]$$
$$+ \sumup_{s\leq 0} C_{iq}^p [\iota_{p,s},
{\cal L}_{j,m}]\epsilon^{q,s-k} +
\sumup_{s\leq 0} C_{iq}^p \iota_{p,s} [\epsilon^{q,s-k}, {\cal L}_{j,m}]$$
$$= -\sumup_{s >0} C_{iq}^p C_{jn}^q \epsilon^{n, s-k-m}\iota_{p,s}
+
\sumup_{s > 0} C_{iq}^p C_{jp}^n \epsilon^{q,s-k}\iota_{n,s+m}$$
$$- \sumup_{s\leq 0} C_{iq}^p C_{jp}^n \iota_{n,
s+m}\epsilon^{q,s-k}
+ \sumup_{s\leq 0} C_{iq}^p C_{jn}^q
\iota_{p,s}\epsilon^{n,s-k-m}$$
$$ =\sumup_{p,q,n} C_{iq}^pC_{jn}^q :\iota_{p,s}\epsilon^{n,s-k-m}: - \sumup_{p,q,n} C_{in}^q C_{jq}^p :\iota_{p,s}\epsilon^{n,s-k-m}: $$
$$- \sumup_{0 < s \leq
m} C_{in}^q C_{jq}^p [\iota_{p,s}, \epsilon^{n,s-k-m}]_+ $$
(Recall $m\geq 0$. Note that the last term vanishes if $m=0$.) Since $[e_{i,k}, e_{j,m}] = [\alpha_{i}, \alpha_j] z^{k+m}$, we have

$${\cal L}([e_{i,k}, e_{j,m}]) = \sumup_q C_{ij}^q {\cal L}_{q, k+m}$$
$$= \sumup_{p,q,s} C_{ij}^q C_{qn}^p :\iota_{p,s} \epsilon^{n, s-k-m}:.$$
In terms of structure constants, the Jacobi identity for $\goth{g}_{\bf C}$
translates into

$$\sumup_{q} (C_{ij}^q C_{qn}^p + C_{jn}^qC_{qi}^p +
C_{ni}^qC_{qj}^p) = 0.$$
Then

$$[{\cal L}_{i,k}, {\cal L}_{j,m}]- {\cal L}([e_{i,k}, e_{j,m}])=
\sumup_{p,q,n} C_{iq}^pC_{jn}^q :\iota_{p,s}\epsilon^{n,s-k-m}:$$
$$- \sumup_{p,q,n} C_{in}^q C_{jq}^p :\iota_{p,s}\epsilon^{n,s-k-m}: -
\sumup_{0 < s \leq m} C_{in}^q C_{jq}^p [\iota_{p,s},
\epsilon^{n,s-k-m}]_+ $$
$$ - \sumup_{p,q,n} C_{ij}^q C_{qn}^p \iota_{p,s} \epsilon^{n, s-k-m}$$
$$= -\sumup_{0 < s \leq m} C_{in}^q C_{jq}^p [\iota_{p,s},
\epsilon^{n,s-k-m}]_+.$$
This shows that, unless $p=n$ and $k= -m$,
$$[{\cal L}_{i.k}, {\cal L}_{j,m}]- {\cal L}([e_{i,k}, e_{j,m}]) = 0,$$
and, when $k=-m$, $p=n$,  we have

$$[{\cal L}_{i,k}, {\cal L}_{j,-k}] - {\cal L}([e_{i,k}, e_{j,-k}]) =
k\cdot\sumup_{q,n} C_{in}^qC_{jq}^n.$$
But, $\sumup_{q,n} C_{in}^qC_{jq}^n = 2c\langle \alpha_i, \alpha_j\rangle$.
The identity obviously does not depend on the assumption that $m>0$. 

Since we had chosen an orthonormal basis with respect to $\langle , \rangle$, $\sumup_{q,n} C_{in}^qC_{jq}^n$ is only non-zero when $i=j$, and 

$${\cal L}([e_{i,k}, e_{i, -k}]) = 0.$$
Thus

$$[{\cal L}_{i,k}, {\cal L}_{i, -k}] = 2c\cdot k.$$
This concludes the proof of Lemma~\ref{hunting}.

\begin{Remark}
\label{remark}
Note that the projective representation of $\goth{A}$, when restricted to
$\goth{g}$, becomes a genuine representation which determines an
action of $G$. There is also a natural rotation action of the circle
$\bf T$ on loops which defines an action of ${\bf
T}\times G$ on the semi infinite forms.
\end{Remark}

Define ${d}: \wedge_\infty^* \rightarrow \wedge_\infty^*$ which
increases degree by $1$ by

\begin{equation}
\label{motherly}
{d} = {1\over 2} \sumup_{i,k}{\cal L}_{i,k} \epsilon^{i,k}.
\end{equation}
Although $d$ is expressed as a sum over all integers $k$, as an
operator on any element of $\wedge^*_{\infty}(\goth{A})$, only finite
number of its terms will be non zero, because, for any element $\omega$ in
$\wedge^*_{\infty}(\goth{A})$, there is an integer $N$ such that
$\omega$ will be annihilated by $\epsilon^{i,k}$ for $k<N$. Also, note that because of our choice of basis ${\cal L}_{i,k}$ and  $\epsilon^{i,k}$ commutes, hence it doesn't matter in which order we write it. Consider the 
$twisted$ operator

$${\tilde d} = {1\over 2} \sumup_{i,k} s_k {\cal L}_{i,k}\epsilon^{i,k},$$
where $s_k= 1$ when $k>0$ and $s_k=-1$ when $k\leq 0$. Take the
adjoint ${\tilde d}^*$ of ${\tilde d}$ and let ${\tilde L}= {\tilde
d}^*{d} + {d}{\tilde d}^*$. Define

$$\Omega = e^{n, 0}\wedge e^{n-1, 0}\wedge \ldots e^{1, 0}\wedge
e^{n,-1}\wedge\ldots$$
(recall that $n$ is the dimension of $\goth{g}_{\bf C}$).
Call this  {\it the vacuum vector} of $\wedge^*_\infty(\goth{A})$.
The Lie algebra cochains $A^*({\goth{a}})$ can be identified with a
subspace of $\wedge_\infty^*(\goth{A})$ by the map $a\mapsto
\epsilon(a)\Omega$. The main statement of this section is the following.

\begin{Proposition}
\label{main-sec2}
The operator $\tilde L$ restricts to the ordinary Laplacian $L= d^*d +
dd^*$ on $A^*(\goth{a})$.
\end{Proposition}
{\bf Proof:} 
The proposition follows from the
following two lemmas.

\begin{Lemma}
\label{jetsdiff}
The operator $d$ restricts to the ordinary Lie algebra
differential (which we will also denote $d$) on the subspace
$A^*(\goth{a})\subset
\wedge_\infty^*(\goth{A})$.
\end{Lemma}

\begin{Lemma}
\label{jetsaddiff}
The adjoint ${\tilde d}^*$ of $\tilde d$ restricts to the adjoint
$d^*$ of the ordinary Lie algebra differential $d$ on
$A^*(\goth{a})\subset \wedge_\infty^*(\goth{A})$.
\end{Lemma}
{\bf Proof of Lemma~\ref{jetsdiff}:}
We need only prove two things. First we will prove that

\begin{equation}
\label{man-eater}
d(\alpha\wedge\omega) = d(\alpha)\wedge\omega \pm
\alpha\wedge d(\omega)
\end{equation}
for $\alpha\in \wedge^*(\goth{a}^*)$ and $\omega\in
\wedge^*_{\infty}(\goth{A})$, where the sign depends on the
degree of $\alpha$. Then we will prove that

\begin{equation}
\label{woman-eater}
d\Omega = 0.
\end{equation}
This will give us Lemma~\ref{jetsdiff} since it will prove that,
for $\alpha\in \wedge^*(\goth{a})$,

$$d(\epsilon(\alpha)\Omega) = \epsilon(d(\alpha))\Omega. $$
Identity~\ref{man-eater} follows, with a bit of calculation, from the
fact that $\epsilon^{i,k}$ anti commutes with any other
$\epsilon^{j,m}$ and the fact that $[\epsilon^{j,m} , {\cal L}_{i,k}]
= \sumup_q C_{iq}^j \epsilon^{q, m-k}$. On the other hand,

$$d\Omega = {1\over 2} \sumup_{i,k} {\cal L}_{i,k}\epsilon^{i,k}\Omega.$$
By the definition of $\Omega$, the only possible non-zero terms are
the ones for which $k>0$. Recall that ${\cal L}_{i,k}$ and
$\epsilon^{i,k}$ commutes, so, all we need to show is that ${\cal
L}_{i,k}\Omega = 0$ for $k>0$. But,

$${\cal L}_{i,k} = \sumup_{p,q,s} C_{iq}^p:\iota_{p,s}\epsilon^{q,s-k}:$$
$$=\sumup_{s\leq 0,p,q}  C_{iq}^p\iota_{p,s}\epsilon^{q,s-k}$$
$$- \sumup_{s>0,p,q} C_{iq}^p\epsilon^{q,s-k}\iota_{p,s}$$
is zero on $\Omega$ since $\epsilon^{q,s-k}$ is zero on $\Omega$ if
$s\leq 0$ and $\iota_{p,s}$ is zero on $\Omega$ if $s>0$.
This concludes the proof of Lemma~\ref{jetsdiff}.

\noindent{\bf Proof of Lemma~\ref{jetsaddiff}:}
For $c_1\in A^*(\goth{a})$ and $c_2\in \wedge_{\infty}^*$,

$$ \langle {\tilde d}^*(c_1\wedge \Omega) , c_2\rangle = \langle
c_1\wedge \Omega, {\tilde d}(c_2)\rangle$$
$$=\langle c_1\wedge\Omega, -{1\over 2}\sumup_{k\leq 0}{\cal
L}_{i,k}\epsilon^{i,k} (c_2) + {1\over 2}\sumup_{k>0}
{\cal L}_{i,k}\epsilon^{i,k}(c_2)\rangle$$
Since all other terms will be killed, we may assume that $c_2$ is a
linear combination of elements of type $c_3\wedge\Omega$ and elements
of type $c_4\wedge \Omega_{j,m}$ where $\Omega_{j,m}$ is the vacuum
vector $\Omega$ with $e^{j,m}$ missing and $c_3, c_4\in
\wedge^*(\goth{a}^*)$. It is enough to show that

\begin{equation}
\label{melgibson}
\langle {\tilde d}^*(c_1\wedge \Omega), c_3\wedge\Omega\rangle = \langle
{d}^*(c_1)\wedge\Omega, c_3\wedge \Omega\rangle
\end{equation}
and that

$$\langle {\tilde d}^*(c_1\wedge \Omega), c_4\wedge\Omega_{j,m}\rangle
= 0.$$ 
For $k\leq 0$, ${\cal L}_{i,k}\epsilon^{i,k}(c_3\wedge\Omega)=0$ and

$$ {1\over 2}\sumup_{k>0} {\cal L}_{i,k}\epsilon^{i,k}(c_3\wedge
\Omega) = {1\over 2}\sumup_{k>0} C_{iq}^p :\iota_{p, s} \epsilon^{q,
s-k}:\epsilon^{i,k} (c_3\wedge \Omega)$$
The sum is over all $p,q,s,i$ as well as $k>0$. Note that, since
$k>0$, the terms, for which $s\leq 0$, are zero (the operator
$\epsilon^{q,s-k}$ is zero on $\Omega$). So,

$$ {1\over 2}\sumup_{k>0} {\cal L}_{i,k}\epsilon^{i,k}(c_3\wedge
\Omega) = d(c_3\wedge\Omega)$$
Hence,

$$\langle {\tilde d}^*(c_1\wedge \Omega), c_3\wedge\Omega\rangle =
\langle {d}^*(c_1\wedge\Omega), c_3\wedge \Omega\rangle = \langle
{d}^*(c_1)\wedge\Omega, c_3\wedge \Omega\rangle$$
(the last equality follows from identity~\ref{man-eater} and
identy~\ref{woman-eater}), which proves the identity~\ref{melgibson}.
On the other hand, 

$${\tilde d}(c_4\wedge \Omega_{j,m}) $$
$$= -{1\over 2}\sumup_{k\leq 0}
{\cal L}_{i,k} \epsilon^{i,k} (c_4\wedge \Omega_{j,m}) + {1\over 2}
\sumup_{k>0} {\cal L}_{i,k}\epsilon^{i,k}(c_4\wedge \Omega_{j,m})$$
$$= -{1\over 2}\sumup_{k\geq 0} {\cal L}_{i,-k}
\epsilon^{i,-k}(c_4\wedge \Omega_{j,m}) + {1\over 2}\sumup_{k>0} {\cal
L}_{i,k}\epsilon^{i,k}(c_4\wedge \Omega_{j,m})$$
$$= -{1\over 2}\sumup_{k\geq 0, s>0} C_{iq}^p
\iota_{p,s}\epsilon^{q,s+k}\epsilon^{i,-k}(c_4\wedge \Omega_{j,m}) +
{1\over 2}\sumup_{k>0, s>0} C_{iq}^p
\iota_{p,s}\epsilon^{q,s-k}\epsilon^{i,k}(c_4\wedge \Omega_{j,m}).$$
Note that, if $s\leq 0$ then the only non-zero terms in $\sumup_{k>0} C_{iq}^p
\iota_{p,s}\epsilon^{q,s+k}\epsilon^{i,-k}(c_4\wedge \Omega_{j,m})$ or
in \\ $\sumup_{k>0} C_{iq}^p
\iota_{p,s}\epsilon^{q,s-k}\epsilon^{i,k}(c_4\wedge \Omega_{j,m})$
which are in $A^*(\goth{a})$ are the ones where $p=q$ or $p=i$. Since
we have chosen an orthonormal basis of $\goth{g}$, these terms are
zero. Now,

$$ \sumup_{k>0, s>0}
C_{iq}^p \iota_{p,s}\epsilon^{q,s-k}\epsilon^{i,k}(c_4\wedge \Omega_{j,m})$$
$$=\sumup_{s>0} \sumup_{k>-s}
C_{iq}^p \iota_{p,s}\epsilon^{q,-k}\epsilon^{i,s+k}(c_4\wedge \Omega_{j,m})$$
Therefore,

$${\tilde d}(c_4\wedge \Omega_{j,m}) = \sumup_{s>0} \sumup_{0> k>-s}
C_{iq}^p \iota_{p,s}\epsilon^{q,-k}\epsilon^{i,s+k}(c_4\wedge \Omega_{j,m})$$
But, for $-s< k \leq 0$,
$\iota_{p,s}\epsilon^{q,-k}\epsilon^{i,s+k}(c_4\wedge \Omega_{j,m})$
can not be contained
in $A^*(\goth{a})$, because $c_4\wedge \Omega_{j,m}$ is missing
$e^{i,k}$ for some $k<0$ and
$\iota_{p,s}\epsilon^{q,-k}\epsilon^{i,s+k}$ can not replace this
missing element. Hence,

$$ \langle {\tilde d}^*(c_1\wedge \Omega), c_4\wedge
\Omega_{j,m}\rangle\ = 0,$$
concluding the proof of Lemma~\ref{jetsaddiff}.

The results in this section show that, to calculate the Lie algebra
cohomology of $\goth{a}$, we need only find the kernel of $\tilde L$.
Since the semi-infinite forms are acted on by ${\bf T}\times G$
(Remark~\ref{remark}), we know that the cochains of $\goth{a}$ are
acted on by ${\bf T}\times G$. In fact we will shortly see that
${\cal L}_{i,0}$ for each $i$ commutes with the operator $d$ and,
hence, the action of ${\bf T}\times G$ on the cochains induces an
action on the cohomology. It follows that the cohomology can be
written as a sum of irreducible representations of ${\bf T}\times
G$. The exact nature of the decomposition will follow from the
explicit formula for $\tilde L$ which will be given in the next
section.

\section{The calculation for the operator ${\tilde L}$}
\label{second}

The main aim of this section is to find a convenient expression of $\tilde
L$ which will enable us to calculate its kernel.

\begin{Proposition}
\label{laplace_formula}
The Laplacian $[{d}, {\tilde d}^*]_+$ is given by

$$[{d}, {\tilde d}^*]_+ = -\sumup_{k>0}c\cdot k\epsilon^{i,k}\iota_{i,k} - \sumup_{k<0}c\cdot k \iota_{i,k}\epsilon^{i,k} + {1\over 2}\sumup_i {\cal
L}_{i,0}^2.$$
\end{Proposition}
{\bf Proof* of Proposition~\ref{laplace_formula}:} 
First note that

$$[{d}, {\tilde d}^*]_+ = [{d}, -{1\over
2}\sumup_{i,k} s_k \iota_{i,k}{\cal L}_{i,k}]_+$$
$$= [{d}, -{1\over 2}(\sumup_{k>0} \iota_{i,k}{\cal L}_{i,-k}
-\sumup_{i,k < 0} \iota_{i,k}{\cal L}_{i,-k} -\sumup_i \iota_{i,0}{\cal
 L}_{i,0})]_+$$
This is equal to

$$-{1\over 2}(\sumup_{k>0} [d, \iota_{i,k}]_+{\cal L}_{i,-k}
-\sumup_{k>0} \iota_{i,k}[{d}, {\cal L}_{i,-k}])$$
$$ +{1\over 2}(\sumup_{k< 0} [d,\iota_{i,k}]_+ {\cal L}_{i,-k}
-\sumup_{k< 0} \iota_{i,k}[d, {\cal L}_{i,-k}])$$
\begin{equation}
\label{lap}
+{1\over 2}(\sumup_i
[d,\iota_{i,0}]_+{\cal L}_{i,0}-\sumup_i \iota_{i,0}[d, {\cal L}_{i,0}])
\end{equation}
To resolve this equation, it is just a matter of identifying $[d,
\iota_{i,k}]_+$ and $[d, {\cal L}_{i,-k}]$. Given a Lie group and its
Lie algebra, the infinitesimal action of the Lie algebra on the Lie
group can be expressed as the Lie derivative which can be written as
the anti-commutator of the differential and the interior product with
respect to the vector fields in the Lie algebra. The operator ${\cal
L}_{i,k}$, with the given basis, on the semi-infinite forms can
likewise be expressed as ${\cal L}_{i,k} = [\iota_{i,k}, {d}]_+$. This
easily follows from the definition of $d$ and the anti-commutators and
commutators calculated in the previous section. Substituting this into
$[{\cal L}_{i,k}, {d}]$, and using the Jacobi identity, it is an easy
calculation to see that

$$[{\cal L}_{i,k}, {d}]= [\iota_{i,k}, {d}^2].$$
We claim that
\begin{Lemma}
\label{diffsquared}
The square ${d}^2$ can be expressed as

$$ {d}^2 = \sumup_{k>0,i} 2c\cdot k \epsilon^{i,k}\epsilon^{i,-k}$$
\end{Lemma}
Lemma~\ref{diffsquared} would imply that

$$[{\cal L}_{i,k}, d] = 2c\cdot k \epsilon^{i,-k}.$$
From this, we see that

\begin{equation}
\label{hey}
[d, {\cal L}_{i,-k}] = 2c\cdot k\epsilon^{i,k}.
\end{equation}
Substituting~\ref{hey} and the identity for $[\iota_{i,k}, {d}]_+$
in~\ref{lap}, we have

$$[{d}, {\tilde d}^*]_+$$
$$= -{1\over 2}\sumup_{k>0} ({\cal L}_{i,k}{\cal L}_{i,-k} - 2c\cdot k
\iota_{i,k}\epsilon^{i,k}) $$
$$+ {1\over 2}\sumup_{k< 0}({\cal L}_{i,k}{\cal L}_{i,-k}- 2c\cdot k
\iota_{i,k}\epsilon^{i,k}) $$
$$+ {1\over 2}\sumup_i {\cal L}_{i,0}{\cal L}_{i,0}$$
Using the commutation rules for ${\cal L}_{i,k}, \iota_{i,k}$ and
$\epsilon_{i,k}$, we have

$$-{1\over 2}\sumup_{k>0} ({\cal L}_{i,k}{\cal L}_{i,-k} - k
\iota_{i,k}\epsilon^{i,k})$$
$$= - {1\over 2}\sumup_{k>0} ( 2c\cdot k + {\cal L}_{i,-k}{\cal L}_{i,k} - 2c\cdot k + 2c\cdot k\epsilon^{i,k}\iota_{i,k})$$
$$=  - {1\over 2}\sumup_{k>0}({\cal L}_{i,-k}{\cal L}_{i,k}
 + 2c\cdot k\epsilon^{i,k}\iota_{i,k})$$
Hence,

$$[{d}, {\tilde d}^*]_+ = {1\over 2}\sumup_{k< 0}{\cal L}_{i,k}{\cal
L}_{i,-k} - {1\over 2}\sumup_{k>0}{\cal L}_{i,-k}{\cal L}_{i,k} $$
$$-\sumup_{k>0}c\cdot k\epsilon^{i,k}\iota_{i,k} - \sumup_{k<0}c\cdot k \iota_{i,k}\epsilon^{i,k} + {1\over 2}\sumup_i {\cal
L}_{i,0}{\cal L}_{i,0}$$
$$= -\sumup_{k>0}c\cdot k\epsilon^{i,k}\iota_{i,k} - \sumup_{k<0}c\cdot
k \iota_{i,k}\epsilon^{i,k} + {1\over 2}\sumup_i {\cal L}_{i,0}^2.$$

To complete the proof of Proposition~\ref{laplace_formula}, we need
only prove Lemma~\ref{diffsquared}.

\noindent{\bf Proof of Lemma~\ref{diffsquared}:}
First of all, we will prove that $d^2$ is a homomorphism on the
$\wedge^*(\goth{a})$ module $\wedge^*_{\infty}(\goth{A})$.
To see this, all we need to show is that

$$[d^2, \epsilon(\alpha)] = 0$$
for any $\alpha\in \wedge (\goth{a})$. In fact,

$$ [d, \epsilon(\alpha)]_+ = \mbox{ multiplication by }d(\alpha)$$
so,

$$[d, [d, \epsilon(\alpha)]_+] = \mbox{ multiplication by }d(d(\alpha)).$$
(Note, $d(d(\alpha)) = 0$.)
But

$$[d^2, \epsilon(\alpha)] = {1\over 2}[d, [d,\alpha]_+] = 0$$
Now let $\Omega_{-k}$ be the element in $\bigwedge_{\infty} =
\wedge^*_{\infty}(\goth{A})$ given by

$$\Omega_{-k} = e^{1,-k}\wedge\ldots \wedge e^{n,-k}\wedge
e^{1,-k-1}\wedge\ldots\wedge e^{n,-k-1}\wedge \ldots$$
($n$ is the dimension of $\goth{g}_{\bf C}$)
and let $\Lambda$ denote the exterior algebra $\wedge(\goth{A})$. Then

$$\mbox{$\bigwedge_{\infty}$ } = \bigcup_{k} \Lambda \Omega_{-k},$$
where $\Lambda \Omega_{-k}$ denotes all elements of the form
$\epsilon(\alpha)\Omega_{-k}$ for $\alpha\in \Lambda$.  Let $I_{-k}$
be ideal of $\Lambda$ generated by the elements $e^{a,-m}$ for
$m\geq k$. Then $\Lambda\Omega_{-k}$ is a $\Lambda/ I_{-k}$
module. Hence, $\bigwedge_{\infty}$ is actually a $\hat\Lambda$
module where $\hat\Lambda$ is the direct limit of $\Lambda/I_{-k}$
as $k$ runs over the positive integers.  We will prove the following
lemma.

\begin{Lemma}
If $T$ is an even degree homomorphism of $\Lambda$-modules then $T$ is
multiplication by an element $\alpha\in {\hat \Lambda}$ of even degree.
\end{Lemma}
{\bf Proof:} 
Note that $\Lambda\Omega_{-k}$ consists of all the
elements $\xi$ of $\bigwedge_{\infty}$ such that $\epsilon(\alpha)\xi
= 0$ for all $\alpha\in I_{-k}$. This shows that
$T(\Lambda\Omega_{-k})\subset \Lambda\Omega$, because $T$ is a
homomorphism on a $\Lambda/I_{-k}$ module $\Lambda\Omega_{-k}$ and
therefore, for all $\alpha\in I_{-k}$

$$\epsilon(\alpha)T(\Lambda\Omega_{-k}) =
T(\epsilon(\alpha)\Lambda\Omega_{-k}) = 0.$$
In fact, for each $k$, $T$ acts on $\Lambda\Omega_{-k}$ as
multiplication by an element $\alpha_{-k}\in
\Lambda/I_{-k}$, where

$$T(\Omega_{-k}) = \alpha_{-k}\Omega_{-k}.$$
But,

$$\Omega_{-k} = \omega_{-k}\Omega_{-k-1},$$
where $\omega_{-k} = e^{1,-k}\wedge\ldots\wedge e^{n,-k}$. So,

$$\alpha_{-k}\Omega_{-k}$$
$$ = T(\Omega_{-k} = T(\omega_{-k} \Omega_{-k-1})
=\omega_{-k}T(\Omega_{-k-1}) =\omega_{-k}\alpha_{-k-1}\Omega_{-k-1}$$
$$=\alpha_{-k-1}\omega_{-k}\Omega_{-k-1} = \alpha_{-k-1}\Omega_{-k}.$$
That means that
$\alpha_{-k}=\alpha_{-k-1}$ in $\Lambda/I_{-k}$. So
$\{\alpha_{-k}\}$ defines an element of $\hat\Lambda$, concluding
the proof of the lemma.

\noindent The lemma proves that $d^2$ is multiplication by some
element $\tilde\omega\in \Lambda^2$. But,

$$[\iota_{i,k}, [\iota_{j,m}, d^2]]_+ = [\iota_{i,k}, ({\cal L}_{j,m}
d- d{\cal L}_{j,m}]$$
$$=[\iota_{i,k}, {\cal L}_{j,m}]d + {\cal L}_{j,m}{\cal L}_{i,k} -
{\cal L}_{i,k}{\cal L}_{j,m} + d[\iota_{i,k}, {\cal L}_{j,m}]$$
$$= \sumup_{p} C_{ij}^p {\cal L}_{p, k+m} - [{\cal L}_{i,k}, {\cal L}_{j,m}]$$
$$= -\delta_{ij}\delta_{k, -m} 2c\cdot k.$$
Suppose ${\tilde \omega} = \sumup_{s,t,u,v} f_{s,t,u,v} e^{s,u}
e^{t,v}$. With a bit of calculation we will be able to see that

$$[\iota_{i,k}, [\iota_{j,m}, \epsilon({\tilde \omega})]]_+ = -2 f_{ijkm}.$$
So,

$${\tilde\omega} = \sumup c\cdot k e^{i,k}e^{i,-k}$$
$$= \sumup_{k>0} 2c\cdot k e^{i,k} e^{i,-k}.$$

The formula of Proposition~\ref{laplace_formula} shows that ${\cal
L}_{i,0}$ commutes with the Laplacian. This proves that on an
irreducible representations of $\goth{g}$ the Laplacian acts by a
scalar. Hence, we need only check how it acts on lowest weight
vectors, which bring us to the main theorem of the paper.

\begin{Theorem}
\label{climax}
The twisted Laplacian ${\tilde L} = [d, {\tilde d}^*]_+$
defined above acts on a irreducible representation $V\subset A^q(\goth{a})$
of ${\bf T}\times G$($\bf T$ is the group of
rotations) with lowest weight $\lambda$
by

$$-\langle \rho, \lambda \rangle + {1\over 2}\parallel \lambda
\parallel^2 - {1\over 2} 2c\cdot k$$
where $k$ is the energy of the lowest weight vector, and $\rho$ is the half sum
of the positive roots of $G$.
\end{Theorem}
{\bf Proof:}
Given the identity~\ref{man-eater}, to determine how
$[d, {\tilde d}^*]$ acts on $A^*(\goth{a})$, it is enough to calculate
it on a vector of the form $v=e^{l,k}\wedge \Omega$. We know from
Proposition~\ref{laplace_formula} that

$${\tilde L}(v)$$
$$= {1\over 2}\sumup_i {\cal L}_{i,0}^2(v) - \sumup_{m>0}c\cdot
m\epsilon^{i,m}\iota_{i,m}(v) - \sumup_{m<0}c\cdot m
\iota_{i,m}\epsilon^{i,m}(v).$$
It is obvious that the third term is zero on $v$. The second term is
also easily seen to be

$$- c\cdot k (v),$$
while to calculate the first term, note that

$$ {\cal L}_{i,0} = \sumup_{p,q,s} C_{iq}^p : \iota_{p,s} \epsilon^{q,s}:.$$
Then, since $\sumup_{p,q,s\leq 0} C_{iq}^p \iota_{p,s}
\epsilon^{q,s}$ kills anything of the form $\beta \wedge\Omega$,

$$ {\cal L}_{i,0}^2 (v)$$
$$ = {\cal L}_{i,0}\{-\sumup_{p,q,s> 0} C_{iq}^p
\epsilon^{q,s}\iota_{p,s}(v)\}$$
$$ = {\cal L}_{i,0}\{-\sumup_{q} C_{iq}^{l} e^{q, k}\wedge \Omega\}$$
$$ = \sumup_{q,r} C_{iq}^{l}C_{ir}^q e^{r,k}\wedge\Omega.$$
But,

$$\sumup_{q,r} C_{iq}^{l}C_{ir}^q $$
$$= \sumup_{q, r} C_{ir}^{q} \langle[\alpha_i, \alpha_q], \alpha_{l}\rangle$$
$$= \sumup_{r} \langle[\alpha_i, [\alpha_i, \alpha_r]], \alpha_{l}\rangle$$
$$= \sumup_{r} \langle[\alpha_i, [\alpha_i, \alpha_l]], \alpha_{r}\rangle$$
(the last equality holds because of the Jacobi identity and the fact
that $\langle [x,y], z\rangle = \langle x, [y,z]\rangle$).  That is,
${1\over 2}\sumup{\cal L}_{i,0}^2$ acts as $- {1\over 2}\sumup_i
\alpha_i^2 = (\mbox{Casmir of }\goth{g})$ (note that the minus sign
comes in because we are acting on the dual space). The action of the
Casmir on a lowest weight vector has been already worked out
(e.g. Section 9.4~\cite{loopgroups})and has the form

$$\sumup_i {\cal L}_{i,0}^2(v)$$
$$ = \{-\langle \rho, \lambda \rangle + {1\over 2}\parallel \lambda
\parallel^2\}(v),$$
where $\lambda$ is the lowest weight for the representation and $\rho$
is the half sum of all positive roots of $G$. This concludes the proof
of Theorem~\ref{climax}.

From Theorem~\ref{climax} it immediately follows that the
cohomology of $\goth{a}$ in any degree is finite dimensional. To
see this, first note that the cochain complex $A^*(\goth{a})$ of
$\goth{a}$ can be divided up according to the energy grading, so that
$A^p(k)(\goth{a})$ are the $p$th cochains of energy $k$. The
differential does not change the energy of the cochain, hence,
$A^*(k)(\goth{a})$ is a subcomplex of $A^*(\goth{a})$. The cochain complex
$A^p(k)(\goth{a})=0$ if $p> k$, i.e. the cohomology $H^*(k)(\goth{a})$
of $A^*(k)(\goth{a})$ is finite dimensional. The cohomology
$H^*(k)(\goth{a})$ is the part of $H^*(\goth{a})$ which is of
energy $k$. On the other
hand, the energy level of the cochains of any one degree is
bounded. To prove this, let us write $H^*(k)(\goth{a}; {\bf C})$ as a
sum of irreducible representations $V_\lambda$ of $G$ with lowest
weight $\lambda$. Theorem~\ref{climax} shows us that the Casmir
operator of $G$ acts on $V_\lambda$ by $P(\lambda) = c\cdot k$. 
Notice that we can choose the basis of $\goth{a}$ so that any
weight of $G$ which occurs in $\wedge^p \goth{a}^*$ is a sum of $p$
roots of $G$. This means that only a finite number of irreducible
representations of $G$ occur in the cohomology of degree $p$.  We 
conclude that $H^p(k)(\goth{a}; {\bf C}) = 0$ if $k > 2\mbox{ sup }
P(\lambda)$ where the supremum is taken over the lowest weight
$\lambda$ of irreducible representations that might occur in
$H^p(\goth{a}; {\bf C})$.

The above argument does not show that there is only one copy of any
one irreducible representation in the sum. But this follows from the
results of the next section.

\section{Conclusion}
\label{third}

Based on the results in the previous sections, we will summarise
the conclusion of this paper in the following theorem.

\begin{Theorem}
\label{one-last}
The $p$th degree cohomology of $\goth{a}$ can be written as a direct
sum of irreducible representations of ${\bf T}\times G$,

$$ H^p(\goth{a}) = \oplus_w V_w,$$
where the sum ranges over elements $w$ of length $p$ in the the quotient
$\frac{{\cal W}_{af}}{\cal W}$ of the affine Weyl group ${\cal W}_{af}$ by
the Weyl group of the finite dimensional Lie group $G$. By
length of an element in the quotient we mean the length of
the shortest representative.
\end{Theorem}
{\bf Proof:}
If $\goth{g}_{\bf C}$ is semi-simple, $\goth{g}_{\bf C} =
\goth{g}_1\oplus \ldots, \oplus \goth{g}_p$ for some $p$ and where
each $\goth{g}_i$ is simple. Then

$$ \goth{a} = \goth{a}_1\oplus \ldots \oplus \goth{a}_p$$
where $\goth{a}_i$ is the corresponding Lie subalgebra of $\goth{A}$
associated to $\goth{g}_i$. Then

$$H^*(\goth{a}) = H^*(\goth{a}_1)\otimes \ldots \otimes
H^*(\goth{a}_p).$$
Hence, it is enough to calculate the cohomology when $\goth{g}_{\bf
C}$ is simple. In this case the formula for the Laplacian
given in ~\ref{climax} can be re-written in a much tidier and familiar
form.

Recall Remark~\ref{remark}. There is an action of ${\bf T}\times G$
and any weight of ${\bf T}\times G$ can be described by a triple ${\bf
\lambda} = (n_1, \lambda, 0)$, where $n_1$ is a number coming from the
action of $\bf T$ as rotation and $\lambda$ is a weight of the $G$
action and the last component is the weight of the central extention
corresponding to the projective representation of $\goth{A}$. Recall that $c$
denotes the coxeter number of $\goth{g}$. And, let $\boldsymbol{\rho}$ denote
the weight $(0, \rho, -c)$ where $\rho$ is a half sum of the positive
roots of $G$. Define

$$\langle(n_1, \lambda_1, b_1), (n_2, \lambda_2, b_2)\rangle =
-n_2b_1 -n_1b_2 + \langle\lambda_1, \lambda_2\rangle,$$
where $\langle,\rangle$ on the right-hand side is the inner product
induced by the killing form on $\goth{g}_{\bf C}$. Then the formula
in Theorem~\ref{climax} for the Laplacian translates into

\begin{equation}
\label{creme}
{\tilde L} = {1\over 2}(\mid\mid
\boldsymbol{\lambda}-\boldsymbol{\rho}\mid\mid^2
-\mid\mid\boldsymbol{\rho}\mid\mid^2),
\end{equation}
on a irreducible representation with lowest weight
$\boldsymbol{\lambda} = (k, \lambda, 0)$ where $\lambda$ is a lowest
weight of $G$ and $k$ denotes the energy of the lowest weight vector.
Note that the basis could have been chosen to be orthonormal on the
subspace $A^*(\goth{a})$ of $\wedge_\infty^*$, consisting of Lie
algebra cochains on $\goth{a}$. We want to find the kernel of $\tilde
L$. 

Note that $\lambda$, as the lowest weight of a representation of
$\goth{g}$ on $A^*(\goth{a})$ can be expressed as a sum of roots (not
necessarily all negative or positive) of $\goth{g}_{\bf C}$. Since
$\goth{a}$ is spanned by $e_{i,k}$ where $k>0$, the circle acts
non-trially and the energy $k$ must always be positive. Hence,
$\boldsymbol{\lambda}$ has to be a sum of positive affine roots which
are not roots of $\goth{g}_{\bf C}$. Given this, the rest of the proof for
Theorem~\ref{one-last}follows from Lemma~\ref{one-more} below.
 
\begin{Proposition}
\label{one-more}
Let $\boldsymbol{\lambda}$ be a sum of positive affine roots.  The
expression~\ref{creme} is zero if and only if the positive affine
roots in the sum $\boldsymbol{\lambda}$ are exactly the set of
positive affine roots turned negative by the shortest representative
of a coset in $\frac{{\cal W}_{af}}{\cal W}$.
\end{Proposition}
{\bf Proof:} 
First assume that equation~\ref{creme} is zero.  Let $\cal P$ denote
the positive alcove. We can choose $w$ in the affine Weyl group ${\cal
W}_{af}$ so that $w(\boldsymbol{\rho}- \boldsymbol{\lambda})\in {\cal
P}$.  Since $\boldsymbol{\lambda}$ can be written as a sum of positive
affine roots, $\boldsymbol{\rho}-w(\boldsymbol{\rho}-
\boldsymbol{\lambda})$ is also a sum of positive affine roots, i.e.,
$\boldsymbol{\rho}-w(\boldsymbol{\lambda}- \boldsymbol{\rho})$ is
positive or zero on anything in the positive alcove. Since
$\boldsymbol{\rho}\in {\cal P}$ and
$w(\boldsymbol{\rho}-\boldsymbol{\lambda})\in {\cal P}$,
$\boldsymbol{\rho} + w(\boldsymbol{\rho}-\boldsymbol{\lambda})\in
{\cal P}$. In fact, because $\boldsymbol{\rho}$ is in the interior of
the positive alcove $\cal P$, so is $\boldsymbol{\rho} +
w(\boldsymbol{\rho}-\boldsymbol{\lambda})$. Hence,

$$ \parallel \boldsymbol{\rho}\parallel^2 - \parallel
\boldsymbol{\rho}-\boldsymbol{\lambda}\parallel^2 $$
$$= \langle \boldsymbol{\rho} - w(\boldsymbol{\rho}-
\boldsymbol{\lambda}), \boldsymbol{\rho} + w(\boldsymbol{\rho} -
\boldsymbol{\lambda})\rangle \geq 0$$ This is only equal to zero if
and only if $\boldsymbol{\rho} - w(\boldsymbol{\rho}-
\boldsymbol{\lambda}) = 0$, and that happens if $\boldsymbol{\lambda}
= \boldsymbol{\rho}-w^{-1}\boldsymbol{\rho}$. Now, it is already known
that $\boldsymbol{\rho}-w^{-1}\boldsymbol{\rho}$ is the sum
$s(\boldsymbol{\lambda})$ of all the positive affine roots which
become negative under the action of $w^{-1}$(see~\cite{loopgroups}
p280).  Furthermore, no other sum of positive affine roots can equal
$\boldsymbol{\lambda}$. Suppose such a sum $\boldsymbol{\alpha}_1 +
\ldots + \boldsymbol{\alpha}_k$ existed. Then
$w^{-1}(s(\boldsymbol{\lambda})) = w^{-1}(\boldsymbol{\alpha}_1) +
\ldots + w^{-1}(\boldsymbol{\alpha}_k)$. The
$w^{-1}(s({\hat\lambda}))$ is a sum of negative roots by construction,
but some of of the $w^{-1}(\boldsymbol{\alpha}_i)$ whould be positive
roots. Those which have turned negative are in the sum
$w^{-1}(s(\boldsymbol{\lambda}))$ and can be
canceled from each side, so we will be left with an identity for which
the left-hand side is a sum of negative roots and the right-hand side
is a sum of positive roots. This is not possible.  Hence, the positive
affine roots in $\boldsymbol{\lambda}$ have to be exactly the set of
positive affine roots which turn negative by some $w^{-1}$ in ${\cal
W}_{aff}$. Recall that $\boldsymbol{\lambda}$ is a sum of positive
affine roots which are not roots of $G$ (see
Remark~\ref{second_remark}. Hence, $w^{-1}$ has to be a representative
of a coset in $\frac{{\cal W}_{aff}}{cal W}$. Since any other element
of ${\cal W}_{aff}$ which belong to the same coset would turn roots of
$G$ negative, $w^{-1}$ is the shortest representative.

On the other hand, assume that $\boldsymbol{\lambda} =
\boldsymbol{\alpha}_1 + \ldots + \boldsymbol{\alpha}_k$ for positive
affine roots $\boldsymbol{\alpha}_i$ and that there exists an element
$w\in {\cal W}$ such that $\boldsymbol{\alpha}_1, \ldots,
\boldsymbol{\alpha}_k$ are exactly the positive roots which become
negative by $w$, then

$${\boldsymbol{\rho}}-\boldsymbol{\lambda} = w^{-1}(\boldsymbol{\rho}),$$
 i.e.,
$$\parallel \boldsymbol{\rho} - \boldsymbol{\lambda} \parallel^2 =
\parallel w^{-1}(\boldsymbol{\rho}) \parallel^2 = \parallel
\boldsymbol{\rho} \parallel^2.$$ In other words,
$P(\boldsymbol{\lambda}) = 0$, concluding the proof of
Proposition~\ref{one-more}.

\section{Comment}
\label{fourth}

Given a finite dimensional Lie group $G$, the loop group of $G$ is the
infinite dimensional Lie group of maps from the circle to $G$. The Lie
algebra of the loop group, which we will call the loop algebra, is a
vector space of maps from the circle to the Lie algebra $\goth{g}$ of
$G$. In this section we wish to make a short comment on the
relationship between the Lie algebra $\goth{a}$ and the loop algebra.
Choose a base point $0$ on the circle $S^1$. Note first that given a
loop in $\goth{g}$ which vanishes at $0$,
we can associate to it its Taylor series at $0$
which could be represented by an infinite formal series

$$ a_1 t + {1\over 2}a_2 t^2 + \ldots$$
where $a_i$ represents the $i$th derivative of the loop at $0$. Denote
the vector space of formal series $\goth{J}$. We will refer to it as
the jet algebra for obvious reasons. It is a Lie algebra. The Taylor
series map is an injective map on the quotient of the loop algebra by
the subalgebra consisting of loops whose derivative vanish to infinite
order at $0$. It is not trivial but a known fact that this map is also
surjective (see~\cite{treves} p390 thm 38.1). Let $\goth{J}_N$ denote
the vector space of polynomials

$$ a_0 + a_1 t + {1\over 2}a_2 t^2 + \ldots + {1\over N!} a_N,$$
with product topology. Take the topology of $\goth{J}$ to be the
inverse limit topology induced by the topology on $\goth{J}_N$. We can
show that the Taylor series map takes the quotient Lie algebra
isomorphically, as topological vector spaces, to $\goth{J}$. Now take
the complexification $\goth{J}_{\bf C}$ of $\goth{J}$. There is a map

$$\goth{a}\mapright{\psi} \goth{J}_{\bf C}$$
which induces a map

$$H^*(\goth{J}_{\bf C}; {\bf C})\mapright{\psi^*} H^*(\goth{a}; {\bf C})$$
in cohomology. Since $\goth{J}$ is inverse limit of $\goth{J}_N$, this
map is injective. Each Lie algebra cochain has an energy level. This
energy level is represented by the sum $\sum k_i$ in the formula for
the Laplacian in the last section. This energy level is not changed by
the differential, and therefore, we may compare the two cohomologies
above on each energy level. Restricted to any one energy level, $\psi$ is
an isomorphism. The formula for the Laplacian again shows that,
for any one cohomology degree, only finite number of energy levels are
involved. Hence, $\psi$ is an isomorphism.

\end{document}